\input amstex
\input amsppt.sty

\def\id{\operatorname{id}}
\def\phi{\varphi}
\def\Re{\operatorname{Re}}
\def\Aut{\operatorname{Aut}}
\def\epsilon{\varepsilon}

\document
\topmatter
\title
Description of all complex geodesics in the symmetrized bidisc
\endtitle
\author
Peter Pflug \& W\l odzimierz Zwonek
\endauthor
\address
Carl von Ossietzky Universit\"at Oldenburg, Institut f\"ur
Mathematik, Postfach 2503, D-26111 Oldenburg, Germany
\endaddress
\email pflug\@mathematik.uni-oldenburg.de
\endemail
\address
Uniwersytet Jagiello\'nski, Instytut Matematyki, Reymonta 4,
30-059 Krak\'ow, Poland
\endaddress
\email zwonek\@im.uj.ed.pl
\endemail
\abstract In the paper we find effective formulas for the complex
geodesics in the symmetrized bidisc.
\endabstract
\thanks The research started while the second author was visiting
Universit\"at Oldenburg in the summer 2003 (supported by DAAD) and
was finished while the stay in November 2003 (supported by DFG
Project: PF 227/8-1). The second author was also supported by the
KBN Grant No  5 P03A 033 21.
\endthanks
\endtopmatter
\document
\subheading{1. Introduction} In the paper we deal with the
problems arising while the study of the so-called symmetrized
bidisc. The study of this domain was initiated by J. Agler, F. B.
Yeh and N. J. Young (see e.g. \cite{Agl-You~1},
\cite{Agl-Yeh-You}). Their interest came from the spectral theory.
However, they found out that some of the results obtained for the
symmetrized bidisc are interesting from the viewpoint of the
geometric function theory, especially from the theory of
holomorphically invariant functions. In particular, the
symmetrized bidisc is an example of a bounded pseudoconvex domain,
which is not biholomorphic to any convex domain and on which the
Lempert function and the Carath\'eodory distance coincide (see
\cite{Agl-You~2}, \cite{Cos~2}). Therefore, it is the first
example showing that the Lempert theorem (i.e. the coincidence of
all holomorphically invariant functions) may extend non-trivially
onto a larger class of domains than that of convex domains. Note
that even more is known, namely, the symmetrized bidisc is not the
exhaustion of domains biholomorphic to convex domains -- see
\cite{Edi~2}.

As a consequence of the above mentioned equality a natural notion
of a complex geodesic arise (see e.g. \cite{Jar-Pfl}, the
definition follows). The main aim of the paper is to give a
description of all complex geodesics of the symmetrized bidisc
(see Theorem 2). As a corollary we get an interesting phenomenon
of the symmetrized bidisc; namely, the complex geodesics in the
symmetrized bidisc extend holomorphically through the boundary
and, unlike in all known cases, their boundaries lie in a very
thin part of the boundary of the symmetrized bidisc (see Corollary
3).

Let us note that the description of complex geodesics together
with the proof of their uniqueness was announced in
\cite{Agl-You~2}.

\bigskip

\subheading{2. Basic notions and definitions in the theory of
holomorphically invariant functions} Let us start with basic
definitions of the holomorphically invariant functions. For more
references on this theory we refer the interested reader to
consult \cite{Jar-Pfl} and \cite{Kob}.

Let $\Bbb D$ denote the unit disk in $\Bbb C$ and let $p$ denote
the Poincar\'e distance on $\Bbb D$. For a domain $D\subset\Bbb
C^n$ and for points $w,z\in D$ we define $$ \gather
c_D(w,z):=\sup\{p(f(w),f(z)):\;f\in\Cal O(D,\Bbb D)\};\\ \tilde
k_D(w,z):=\inf\{p(\lambda_1,\lambda_2):\; \text{ there is
$f\in\Cal O(\Bbb D,D)$ with $f(\lambda_1)=w$, $f(\lambda_2)=z$}
\};\\ k_D:=\text{the largest pseudodistance not larger than
$\tilde k_D$}.
\endgather
$$ We call $c_D$ (respectively, $k_D$) the {\it Carath\'eodory }
(respectively, {\it Kobayashi}) {\it pseudodistance}. $\tilde k_D$
is called the {\it Lempert function}. The following simple
inequalities are crucial: $$ c_D\leq k_D\leq \tilde k_D.$$

The important theorem of Lempert states that if $D$ is
biholomorphically equivalent to a convex domain then the
inequalities above become equalities (see \cite{Lem}), i.e., $$
\text{for any domain $D$ biholomorphic to a convex domain we have
} c_D=\tilde k_D. $$

A mapping $\phi\in\Cal O(\Bbb D,D)$ is called a {\it $\tilde
k_D$-geodesic } for $(w,z)$, $w\neq z$, if $\phi(\lambda_1)=w$,
$\phi(\lambda_2)=z$ and $p(\lambda_1,\lambda_2)=\tilde k_D(w,z)$
for some $\lambda_1,\lambda_2\in \Bbb D$. In other words, $\tilde
k_D$-geodesics are those mappings for which the infimum in the
definition of $\tilde k_D$ is attained.

If $D$ is a taut domain (the domain $D$ is {\it taut } if for any
sequence $(\phi_{\nu})_{\nu=1}^{\infty}\subset\Cal O(\Bbb D,D)$
there is a subsequence converging locally uniformly on $\Bbb D$ to
$\phi\in\Cal O(\Bbb D,D)$ or there is a subsequence
$(\phi_{\nu_k})_{k=1}^{\infty}$ such that for any compact
$L\subset D$ and for any compact $K\subset\Bbb D$ there is a
$k_0\in\Bbb N$ such that for any $k\in\Bbb N$, $k\geq k_0$
$\phi_{\nu_k}(K)\cap L=\emptyset$) then for any $w\neq z$, $w,z\in
D$ there is a $\tilde k_D$-geodesic for $(w,z)$. It is simple to
conclude from the Schwarz-Pick lemma that if $\phi\in\Cal O(\Bbb
D,D)$, $w,z\in D$, $w\neq z$, $\phi(\lambda_1^0)=w$,
$\phi(\lambda_2^0)=z$ and $p(\lambda_1^0,\lambda_2^0)=c_D(w,z)$
for some $\lambda_1^0,\lambda_2^0\in \Bbb D$, then
$c_D(\phi(\lambda_1),\phi(\lambda_2))=\tilde
k_D(\phi(\lambda_1),\phi(\lambda_2))= p(\lambda_1,\lambda_2)$ for
any $\lambda_1,\lambda_2\in \Bbb D$. This observation leads us to
the following definition. A mapping $\phi\in \Cal O(\Bbb D,D)$ is
called {\it a complex geodesic} (in $D$) if $$
c_D(\phi(\lambda_1),\phi(\lambda_2))=p(\lambda_1,\lambda_2), $$
for any $\lambda_1,\lambda_2\in \Bbb D$. Certainly, if $\phi$ is a
complex geodesic then $c_D(\phi(\lambda_1),\phi(\lambda_2))=\tilde
k_D(\phi(\lambda_1),\phi(\lambda_2))$, $\lambda_1,\lambda_2\in\Bbb
D$.

The problem of finding explicit formulas for complex geodesics (or
$\tilde k_D$-geodesics) is, in general, very difficult. One of
very few non-trivial examples for which the formulas for complex
geodesics are known completely are convex complex ellipsoids (see
\cite{Jar-Pfl-Zei}) or more generally convex generalized
pseudoellipsoids (see \cite{Zwo}). Another example of a convex
domain, where the complex geodesics are known is the so-called
minimal ball (see \cite{Pfl-Yous}). Without the assumption of
convexity only necessary forms of $\tilde k_D$-geodesics are known
(see \cite{Edi~1} also \cite{Pfl-Zwo}).

\bigskip

\subheading{3. Basic properties of the symmetrized bidisc and
statement of results} Now let us start the study of the
symmetrized bidisc. Define $$ \pi:\Bbb
D\owns(\lambda_1,\lambda_2)\mapsto(\lambda_1+\lambda_2,\lambda_1\lambda_2)\in\Bbb
C^2. $$

Put $G_2:=\pi(\Bbb D^2)$. It is easy to see that $G_2$ is a
domain. We call the domain $G_2$ the {\it symmetrized bidisc}. It
is easy to calculate that $$ G_2=\left\{(s,p)\in\Bbb
C^2:\max\left\{\left|\frac{s\pm\sqrt{s^2-4p}}{2}\right|\right\}<1\right\}.\tag{1}
$$ It follows from this description that $G_2$ is a bounded
hyperconvex domain in $\Bbb C^2$ (we recall that the bounded
domain $D\subset\Bbb C^n$ is {\it hyperconvex} if there is a
continuous plurisubharmonic function $u:D\mapsto(-\infty,0)$ such
that $\lim\sb{D\owns z\to\partial D}u(z)=0$). In particular, $G_2$
is taut, so for any pair of points $w,z\in D$, $w\neq z$ there is
a $\tilde k_{G_2}$-geodesic for $(w,z)$. Let us recall that there
is a number of other possible descriptions of the symmetrized
bidisc (see e.g. \cite{Agl-You~1}, \cite{Agl-You~2} and
\cite{Agl-Yeh-You}). Nevertheless, since we make use only of the
above mentioned descriptions, therefore we do not recall them
here.

Note that $\pi:\Bbb D^2\mapsto G_2$ is a proper holomorphic
mapping with multiplicity $2$ and $\pi:\Bbb
D^2\setminus\triangle\mapsto G_2\setminus S$ is a holomorphic
covering, $\triangle:=\{(\lambda,\lambda):\lambda\in\Bbb D\}$,
$S:=\{(2\lambda,\lambda^2):\lambda\in\Bbb D\}$.

Our first aim is a description of $k_{\tilde G_2}$-geodesics for
the points $((0,0),(s_0,p_0))$ with $(s_0,p_0)\neq(0,0)$. The
result below may be treated as the generalization of the 'Schwarz
Lemma' for the symmetrized bidisc as formulated in
\cite{Agl-You~1} (we give an alternative description of geodesics
and the uniqueness of the solution of the extremal problem -- i.e.
the uniqueness of geodesics passing through $(0,0)$). Moreover,
the description of the $\tilde k_{G_2}$-geodesics seems to us to
be simpler and more natural than that in \cite{Agl-You~1}.

\proclaim{Theorem 1} Let $\phi:\Bbb D\mapsto G_2$ be a holomorphic mapping
such that $\phi(0)=(0,0)$, $\phi(\sigma^2)=(s_0,p_0)\neq(0,0)$ for some
$\sigma\in\Bbb D$. Then the following are equivalent: $$ \gather \phi \text{
is a $\tilde k_{G_2}$-geodesic for $((0,0),(s_0,p_0))$};\tag{2}\\
\text{there is a Blaschke product of degree less than or equal to two}
\tag{3}\\ \text{ with $B(0)=0$ and } \phi(\lambda)=
(B(\sqrt{\lambda})+B(-\sqrt{\lambda}),B(\sqrt{\lambda})B(-\sqrt{\lambda})),\;\lambda\in\Bbb
D.\endgather $$ Moreover, any mapping satisfying \thetag{2} (or \thetag{3})
is a complex geodesic.

Additionally, the $\tilde k_{G_2}$-geodesics for
$((0,0),(s_0,p_0))$ are unique (up to automorphisms of $\Bbb D$).
\endproclaim

\subheading{Remark} It is easy to see that the function $\phi$
defined in \thetag{3} is a holomorphic function on $\Bbb D$,
$\phi(0)=(0,0)$ and its image lies in $G_2$. Moreover, we may
write down explicitly all the possible forms of $\phi$ satisfying
\thetag{3}. Namely, if $B(\lambda)=\tau\lambda$, $\lambda\in\Bbb
D$, $|\tau|=1$ then $$ \phi(\lambda)=(0,-\tau^2\lambda),\;
\lambda\in\Bbb D.\tag{4} $$ If
$B(\lambda):=\tau\lambda\frac{\lambda-\alpha}{1-\bar\alpha\lambda}$,
$\lambda\in\Bbb D$, $\tau\in\partial\Bbb D$, $\alpha\in\Bbb D$
then $$
\phi(\lambda)=\left(\frac{2\tau\lambda(1-|\alpha|^2)}{1-\bar\alpha^2\lambda},
\tau^2\lambda\frac{\lambda-\alpha^2}{1-\bar\alpha^2\lambda}\right),\;\lambda\in\Bbb
D.\tag{5} $$

As already mentioned, the Lempert function and the Carath\'eodory
distance of the symmetrized bidisc coincide (see \cite{Agl-You~2},
\cite{Cos~2}), i.e. $$c_{G_2}=\tilde k_{G_2}.$$ This, together
with the fact that $G_2$ is not biholomorphic to any convex domain
(see \cite{Cos~1}), makes the domain interesting from the
viewpoint of the theory of holomorphically invariant functions. In
the theorem below we shall make use of the above mentioned
equality to give a description of all $\tilde k_{G_2}$-geodesics
(or equivalently, complex geodesics) in $G_2$.

Before we formulate the result let us recall how we can define
automorphisms of $G_2$, which map $(0,0)$ onto an arbitrary point
of $S$ (see \cite{Agl-You~1}). Namely for the fixed $\alpha\in\Bbb
D$ we define the automorphism of $G_2$ (mapping $(0,0)$ onto
$(2\alpha,\alpha^2)$ as follows. Let
$b_{\alpha}(\lambda):=\frac{\alpha-\lambda}{1-\bar\alpha\lambda}$,
$\lambda\in\Bbb D$. Assume that
$(s,p)=\pi(\lambda_1,\lambda_2)=\pi(\lambda_2,\lambda_1)$, where
$(\lambda_1,\lambda_2)\in\Bbb D^2$. Define the automorphism
$B_{\alpha}$ of $G_2$ as follows: $$
B_{\alpha}(s,p):=\pi(b_{\alpha}(\lambda_1),b_{\alpha}(\lambda_2)).
$$

Therefore, in view of Theorem 1, the whole problem with the
description of complex geodesics reduces to the characterization
of geodesics omitting the set $S$ (i.e. such that $\phi(\Bbb
D)\cap S=\emptyset$).

The description of geodesics follows.

\proclaim{Theorem 2} Let $\phi:\Bbb D\mapsto G_2$. Then $\phi$ is
a complex geodesic if and only if one of two possibilities holds

\item{{\rm (i)}} if $\phi(\Bbb D)\cap S\neq\emptyset$, then $\phi=\tilde\phi\circ b_{\alpha}$, where $\alpha\in\Bbb
D$,
$$\tilde\phi(\lambda)=\pi(B(\sqrt{\lambda}),B(-\sqrt{\lambda})),\;\lambda\in\Bbb
D,$$ and $B$ is a Blaschke product of degree one or two;

\item{{\rm (ii)}} if $\phi(\Bbb D)\cap S=\emptyset$, then $\phi=\pi\circ f$, where $f=(f_1,f_2)$, $f_1,f_2\in\Aut\Bbb
D$ and $f_1-f_2$ has no root in $\Bbb D$. \endproclaim

As an immediate corollary of Theorem 2 we get the following
interesting property of complex geodesics of the symmetrized
bidisc.

\proclaim{Corollary 3} Let $\phi:\Bbb D\mapsto G_2$ be a complex
geodesic. Then $\phi$ extends holomorphically onto a neighborhood
of $\bar\Bbb D$ and $\phi(\partial\Bbb D)\subset\pi(\partial\Bbb
D\times\partial\Bbb D)$.
\endproclaim

As a corollary of  Theorem 1 (uniqueness), Theorem 2 and the
existence of automorphisms $B_{\alpha}$ we get the following.

\proclaim{Corollary 4} Let $\phi:\Bbb D\mapsto G_2$ be a complex
geodesic. Then one of the three possibilities occurs:

\noindent{{\rm (i)}} $\phi(\Bbb D)\cap S=\emptyset$;

\noindent{{\rm (ii)}} $\#(\phi(\Bbb D)\cap S)=1$;

\noindent{{\rm (iii)}} $\phi(\Bbb D)=S$.

Moreover, each of the possibilities \thetag{i}--\thetag{iii} is
satisfied by some complex geodesic.
\endproclaim

In the proofs of Theorems 1 and 2 essential role will be played by
special rational functions among which we shall choose candidates
for the extremals in the definition of the Carath\'eodory
distance. Namely, put

$$ \gather F_0(s,p)=p,\; (s,p)\in G_2,\tag{6}\\
F_{\omega}(s,p):=\frac{2p-\omega s}{2-\bar\omega s},\;(s,p)\in
G_2, \text{ where $\omega\in\partial\Bbb D$}.\tag{7}
\endgather
$$ One may verify that $F_{\omega}(G_2)\subset\Bbb D$,
$\omega\in\partial\Bbb D\cup\{0\}$ (see \cite{Agl~You}).

Now we go on to proofs of Theorems 1 and 2.

\bigskip

\subheading{4. Proof of Theorem 1}

\demo{Proof of Theorem 1} First we show that any function $\phi$
defined in \thetag{3} is a complex geodesic. In particular, it
will show the implication ($\thetag{3}\implies\thetag{2}$).

Let us fix $\phi$ as in \thetag{3}. It is sufficient to find $$
F\in\Cal O(G_2,\Bbb D) \text{ such that $F\circ\phi=r$}, $$ where
$r$ is some rotation.

In fact, we shall show that the function $F$ may be chosen as the
one as in \thetag{6} or \thetag{7}.

Let us consider the case when $B(\lambda)=\tau\lambda$,
$\lambda\in\Bbb D$, where $\tau\in\partial\Bbb D$. Then
$\phi(\lambda)=(0,-\tau^2\lambda)$, $\lambda\in\Bbb D$, and then
$F(s,p):=p$, $(s,p)\in D$, satisfies the desired property.

Therefore, we are left with the case
$B(\lambda):=\tau\lambda\frac{\lambda-\alpha}{1-\bar\alpha\lambda}$,
$\lambda\in\Bbb D$, where $\tau\in\partial\Bbb D$ and
$\alpha\in\Bbb D$. In this case the function $F$ will be as in
\thetag{7}.

First we consider the case $\alpha=0$, i.e.
$\phi(\lambda)=(2\tau\lambda,\tau^2\lambda^2)$, $\lambda\in\Bbb
D$. It is sufficient to take $\omega=1$. So assume that
$\alpha\neq 0$. Recall that then $$
\phi(\lambda)=\left(\frac{2\tau\lambda(1-|\alpha|^2)}{1-\bar\alpha^2\lambda},
\tau^2\lambda\frac{\lambda-\alpha^2}{1-\bar\alpha^2\lambda}\right),\;\lambda\in\Bbb
D. $$ Note that in view of the Schwarz Lemma it is sufficient to
choose $\omega$ so that $|F(\phi(\alpha^2))|=|\alpha|^2$ (because
$F(\phi(0))=0$ and $\alpha^2\neq 0$). Let us calculate (substitute
in the last equality
$\omega:=\frac{|\alpha|^2}{\bar\alpha^2}\tau$) $$
|F(\phi(\alpha^2))|=|F(\frac{2\tau\alpha^2}{1+|\alpha|^2},0)|=
\frac{2|\alpha|^2}{(1+|\alpha|^2)\left|2-\bar\omega\frac{2\tau\alpha^2}{1+|\alpha|^2}\right|}=
|\alpha|^2. $$ This completes the proof of the fact that any
function $\phi$ defined by \thetag{3} is a complex geodesic.

\bigskip

Now we show the following:

for any $(s_0,p_0)\in D$, $(s_0,p_0)\neq(0,0)$ there is a mapping
$\phi$ as defined in \thetag{3} such that $\phi$ goes through
$(s_0,p_0)$.

Let $(s_0,p_0)=\pi(t_1,t_2)$. Obviously $(t_1,t_2)\neq(0,0)$. It
is sufficient to show that there is a finite Blaschke product of
degree one or two and $\sigma\in\Bbb D\setminus\{0\}$ such that
$B(0)=0$, $B(\sigma)=t_1$, and $B(-\sigma)=t_2$.

In the case $t_2=-t_1$ it is sufficient to take $B=\id$,
$\sigma=t_1$. So assume that $t_2\neq-t_1$. It is suffices to find
$\tau\in\partial\Bbb D$, $\sigma\in(0,1)$, and $\alpha\in\Bbb D$
such that $$
\tau\frac{\sigma-\alpha}{1-\bar\alpha\sigma}=\frac{t_1}{\sigma},\;\tau\frac{-\sigma-\alpha}{1+\bar\alpha\sigma}=
\frac{t_2}{-\sigma}. $$ Now comparing the values
$p(\sigma,-\sigma)$ and
$p(\frac{t_1}{\sigma},\frac{t_2}{-\sigma})$ as
$\sigma\in(\max\{|t_1|,|t_2|\},1)$ and making use of the Darboux
property we get the existence of a
$\sigma\in(\max\{|t_1|,|t_2|\},1)$ such that
$p(\sigma,-\sigma)=p(\frac{t_1}{\sigma},\frac{t_2}{-\sigma})$.
This easily implies the existence of $\tau$ and $\alpha$ as
desired.

\bigskip

Now we prove the implication ($\thetag{2}\implies\thetag{3}$)
under the additional assumption that $$ \psi^{-1}(0)=\{0\},\text{
where
$\psi(\lambda):=\phi_1^2(\lambda)-4\phi_2(\lambda),\;\lambda\in\Bbb
D$}.\tag{8} $$ Note that $\psi$ is a bounded holomorphic function.
We may write $\psi(\lambda)=\lambda^k\tilde\psi(\lambda)$,
$\lambda\in\Bbb D$, where $\tilde\psi:\Bbb D\mapsto\Bbb C$ is a
bounded holomorphic function and $\tilde\psi(\lambda)\neq 0$,
$\lambda\in\Bbb D$.

Then, in view of the description of $G_2$ as in \thetag{1}. $$
\left|\phi_1(\lambda^2)+\lambda^k\sqrt{\tilde\psi(\lambda^2)}\right|<2,\;\lambda\in\Bbb
D, $$ where the root taken above is chosen arbitrarily so that the
function $B$ defined $$
B(\lambda):=\frac{\phi_1(\lambda^2)+\lambda^k\sqrt{\tilde\psi(\lambda^2)}}{2},\;\lambda\in\Bbb
D, $$ is holomorphic on $\Bbb D$. Note that $B(0)=0$ and $|B|<1$
on $\Bbb D$. Then we may write $$
\lambda^k\sqrt{\tilde\psi(\lambda^2)}=2B(\lambda)-\phi_1(\lambda^2),\lambda\in\Bbb
D. $$ Taking the square in the last equality we get $$
\psi(\lambda^2)=\lambda^{2k}\tilde\psi(\lambda^2)=4B^2(\lambda)-4B(\lambda)\phi_1(\lambda^2)+
\phi_1^2(\lambda^2),\;\lambda\in\Bbb D. $$ Consequently, $$
B(\lambda)\phi_1(\lambda^2)-B^2(\lambda)=\phi_2(\lambda^2),\;\lambda\in\Bbb
D.\tag{9} $$ It follows from \thetag{9} that for any
$\lambda\in\Bbb D$ $$
B(\lambda)\phi_1(\lambda^2)-B^2(\lambda)=\phi_2(\lambda^2)=
B(-\lambda)\phi_1(\lambda^2)-B^2(-\lambda),\;\lambda\in\Bbb D, $$
or $$
(B(\lambda)-B(-\lambda))(\phi_1(\lambda^2)-(B(\lambda)+B(-\lambda)))=0,\;\lambda\in\Bbb
D. $$ In view of the identity principle for holomorphic functions
two cases are possible:

Case (I). $B(\lambda)=B(-\lambda)$, $\lambda\in\Bbb D$.

Case (II). $\phi_1(\lambda^2)=B(\lambda)+B(-\lambda)$,
$\lambda\in\Bbb D$.
\bigskip

In Case (I) we easily get the existence of a function $B_1\in\Cal
O(\Bbb D,\Bbb D)$ such that $B_1(\lambda^2)=B(\lambda)=
\frac{\phi_1(\lambda^2)+\lambda^k\sqrt{\tilde\psi(\lambda^2)}}{2}$,
$\lambda\in\Bbb D$. Taking the other holomorphic branch in the
formula for $B_1$ we may define analogously the function
$B_2\in\Cal(\Bbb D,\Bbb D)$ such that $$
\phi=\pi\circ(B_1,B_2).\tag{10} $$

Then it follows from our earlier considerations that there is some
$F\in\Cal O(G_2,\Bbb D)$ such that $F(0,0)=0$ and
$|F(\phi(\lambda))|=|\lambda|$, where $F$ is as in \thetag{6} or
\thetag{7}. If $F$ is as in \thetag{6} then, in view of the
Schwarz Lemma, $$
|F(\pi(B_1(\lambda),B_2(\lambda)))|=|B_1(\lambda)||B_2(\lambda)|\leq|\lambda|^2,\;\lambda\in\Bbb
D $$ -- contradiction.

In the second case there is some $\tau\in\partial\Bbb D$ such that
$$
\frac{2B_1(\lambda)B_2(\lambda)-\omega(B_1(\lambda)+B_2(\lambda))}{2-\bar\omega(B_1(\lambda)+B_2(\lambda))}
=\tau\lambda,\;\lambda\in\Bbb D;$$ so
$2B_1(\lambda)B_2(\lambda)-\omega(B_1(\lambda)+B_2(\lambda))=2\tau\lambda-\tau\lambda
\bar\omega(B_1(\lambda)+B_2(\lambda))$, $\lambda\in\Bbb D$.
Differentiating at $0$ we get (recall that $B_1(0)=B_2(0)=0$) $$
-\omega(B_1^{\prime}(0)+B_2^{\prime}(0))=2\tau. $$ The last
equality together with the Schwarz inequality applied to $B_1$ and
$B_2$ shows that
$B_1(\lambda)=-\bar\omega\tau\lambda=B_2(\lambda)$,
$\lambda\in\Bbb D$, so
$\phi(\lambda)=(-2\bar\omega\tau\lambda,\bar\omega^2\tau^2\lambda^2)$,
$\lambda\in\Bbb D$, which shows that $\phi$ is as in \thetag{4}
(and contradicts \thetag{8}), which finishes the proof in this
case.

In Case (II) we easily get that $$
\phi_1(\lambda)=B(\sqrt{\lambda})+B(-\sqrt{\lambda}),\;\lambda\in\Bbb
D. $$ Then in view of \thetag{9} we get that $$
\phi_2(\lambda)=B(\sqrt{\lambda})B(-\sqrt{\lambda}),\;\lambda\in\Bbb
D. $$ Therefore, we are left with the proof of the fact that $B$
is a Blaschke product of degree less than or equal to $2$.

Let $(s_0,p_0)=\pi(t_1,t_2)=\pi(t_2,t_1)$, where $t_1,t_2\in\Bbb
D$. Let $\phi(\sigma^2)=(s_0,p_0)$. The assumption that $\phi$ is
a $\tilde k_D$-geodesic for $((0,0),(s_0,p_0))$ leads us to the
following observation: $B(\sigma)=t_1$, $B(-\sigma)=t_2$ (or vice
versa, which, however, may be dealt analogously). Moreover, the
function $B$ is extremal in the following sense: $B\in\Cal O(\Bbb
D,\Bbb D)$, $B(0)=0$, $B(\sigma)=t_1$, $B(\sigma_2)=t_2$  and
there is no holomorphic function $f:\Bbb D\mapsto\Bbb C$ such that
$f(\Bbb D)\subset\subset\Bbb D$ and $f(0)=0$, $f(\sigma)=t_1$,
$f(-\sigma)=t_2$). Therefore, (see e.g. \cite{Edi~1}) $B$ must be
a Blaschke product of degree smaller than or equal to two. So the
implication ($\thetag{2}\implies\thetag{3}$) (under the assumption
\thetag{8}) is finished.

Now we show that the assumption that $\phi$ is a $\tilde
k_{G_2}$-geodesic and the inequality $\psi^{-1}(0)\neq\{0\}$
implies that $\phi(\lambda)=(2\tau\lambda,\tau^2\lambda^2)$,
$\lambda\in\Bbb D$, for some $\tau\in\partial \Bbb D$. In fact,
assume that $\psi(\lambda_0)=0$ for some $\lambda_0\in\Bbb
D\setminus\{0\}$. Consequently,
$\phi_1^2(\lambda_0)=4\phi_2(\lambda_0)$. It follows from our
considerations that $\phi$ must be a complex geodesic and that
there is an $F\in\Cal O(G_2,\Bbb D)$ such that $F\circ\phi$ is a
rotation, where $F$ is as in \thetag{6} or \thetag{7}. In the
first case $|\phi_2(\lambda)|=|\lambda|$, $\lambda\in\Bbb D$. But
because of the Schwarz Lemma we have $|\phi_1^2(\lambda)|\leq
4|\lambda|^2$, which together with the equality
$|\phi_1^2(\lambda_0)|=4|\phi_2(\lambda_0)|=4|\lambda_0|$ gives
the contradiction. So assume that there are
$\omega,\tau\in\partial\Bbb D$ such that $$
\tau\lambda=\frac{2\phi_2(\lambda)-\omega\phi_1(\lambda)}{2-\bar\omega\phi_1(\lambda)}.\tag{12}
$$ Then we get the following equality (substitute
$\lambda=\lambda_0$ in \thetag{12}) $$
\phi_1^2(\lambda_0)+2\phi_1(\lambda_0)(\bar\omega\tau\lambda_0-\omega)-4\tau\lambda_0=0.
$$ Solving the above equality with respect to $\phi_1(\lambda_0)$
(note that $|\phi_1(\lambda_0)|<2$) we get that
$\phi(\lambda_0)=-2\bar\omega\tau\lambda_0$. Applying the Schwarz
Lemma to $\phi_1:\Bbb D\mapsto 2\Bbb D$ we get that
$\phi_1(\lambda)=-2\bar\omega\tau\lambda$, $\lambda\in\Bbb D$.
Substituting this formula in \thetag{12} it follows that
$\phi_2(\lambda)=\bar\omega^2\tau^2\lambda^2$, $\lambda\in\Bbb D$,
which shows that $\phi$ is as in \thetag{4}.
\bigskip

To finish the proof of the theorem it is sufficient to show the
uniqueness. Let $\phi$ be a $\tilde k_{G_2}$-geodesic for
$((0,0),(s_0,p_0))$ and $\phi(\sigma^2)=(s,p)= \pi(t_1,t_2)$. Then
$\phi$ is given by some Blaschke product of degree less than or
equal to two. Moreover, the Blaschke product is such that $B(0)=0$
and, without loss of generality, $B(\sigma)=t_1$,
$B(-\sigma)=t_2$. But the values at three points ($0$, $\sigma$
and $-\sigma$) determine uniquely such Blaschke products. This
completes the proof of the uniqueness and the proof of the whole
theorem. \qed
\enddemo

\subheading{Remark} Note that in the proof of Theorem 1 we did not
make use of the equality $\tilde k_{G_2}=c_{G_2}$. In fact, we
even showed the equality $$\tilde
k_{G_2}((0,0),(s_0,p_0))=c_{G_2}((0,0),(s_0,p_0)) \text{ for any
$(s_0,p_0)\in G_2$}.$$ Assuming the equality $\tilde
k_{G_2}=c_{G_2}$ from the very beginning we may simplify the proof
a little. However, without this result our proof is not much
longer but much more self-contained.

\bigskip

\subheading{5. Proof of Theorem 2} In this part of the paper we
provide the proof of Theorem 2. In contrast to the proof of
Theorem 1 we make use of the equality $\tilde k_{G_2}=c_{G_2}$. In
particular, we know that the notions of $\tilde k_{G_2}$-geodesics
and complex geodesics coincide in the symmetrized bidisc.

\demo{Proof of Theorem 2} Note that the composition
$b_{\beta}\circ B$, where $\beta\in\Bbb D$ and $B$ is a Blaschke
product of degree one or two, is a Blaschke product of degree one
or two. Therefore, in view of Theorem 1 and because of the
existence of automorphisms $B_{\alpha}$ mapping $(0,0)$ to an
arbitrarily chosen point from $S$, to prove Theorem 2 it is
sufficient to study the mappings $\phi$ omitting $S$.

Since we need this fact, note that we already know that if $\phi$
is a complex geodesic in $G_2$ such that $\phi(\Bbb D)\cap
S\neq\emptyset$ then $\phi(\partial\Bbb D)\subset\pi(\partial\Bbb
D\times\partial\Bbb D)$.

To finish the proof of Theorem 2 it is sufficient to show the
following equivalence.

Let $\phi:\Bbb D\mapsto G_2$, $\phi(\Bbb D)\cap S=\emptyset$. Then
$\phi$ is a complex geodesic if and only if $\phi=\pi\circ f$,
where $f=(f_1,f_2)$, $f_1,f_2\in\Aut\Bbb D$ and $f_1-f_2$ has no
root in $\Bbb D$.

($\Rightarrow$): Let $\phi$ be a complex geodesic such that
$\phi(\Bbb D)\cap S=\emptyset$. Then there is a holomorphic
function $f:\Bbb D\mapsto\Bbb D^2$ such that $f(\Bbb D)\cap
\triangle=\emptyset$ and $\phi=\pi\circ f$. One may easily verify
that $f$ is $\tilde k_{\Bbb D^2}$-geodesic, so it is a complex
geodesic in $\Bbb D^2$. It is sufficient to show that
$f=(f_1,f_2)$, where $f_1,f_2\in\Aut\Bbb D$. Since $f$ is a
complex geodesic we easily obtain that at least one of its
components is from $\Aut\Bbb D$. Suppose that the other component
is not from $\Aut\Bbb D$. Without loss of generality let
$f_1\in\Aut\Bbb D$ and $f_2\not\in\Aut\Bbb D$. Fix
$\sigma\in(0,1)$. Then there is a function $\tilde f_2\in\Cal
O(\bar \Bbb D,\Bbb D)$ such that $\tilde f_2(0)=f_2(0)$, $\tilde
f_2(\sigma)=f_2(\sigma)$ and $\tilde f_2(\Bbb D)\subset\subset\Bbb
D$. Then $(f_1,\tilde f_2)$ is a $\tilde k_{\Bbb D^2}$-geodesic
for $(f(0),f(\sigma))$ (and consequently, a complex geodesic) and
$\pi\circ (f_1,\tilde f_2)$ is a complex geodesic in $G_2$. But
$f_1-\tilde f_2$ has one root in $\Bbb D$ (use Hurwitz theorem),
so $\tilde\phi:=\pi\circ(f_1,\tilde f_2)$ intersects $S$, which
however contradicts the description of complex geodesics passing
through $S$ because $\tilde\phi(\partial\Bbb
D)\not\subset\pi(\partial\Bbb D\times\partial\Bbb D)$.

\bigskip

($\Leftarrow$): Let $\phi=\pi\circ(f_1,f_2)$ be such that
$f_1,f_2\in\Aut\Bbb D$ and the equation $$
f_1(\lambda)=f_2(\lambda)\tag{13} $$ has no root for
$\lambda\in\Bbb D$. Without loss of generality we may assume that
$f_1(\lambda)=\lambda$, $\lambda\in\Bbb D$ and
$f_2(\lambda)=\tau\frac{\lambda-\alpha}{1-\bar\alpha\lambda}$,
$\lambda\in\Bbb D$, where $\tau\in\partial\Bbb D$ and
$\alpha\in\Bbb D\setminus\{0\}$. One may easily verify that the
nonexistence of roots of the equation \thetag{13} is equivalent to
the inequality $$ |1-\tau|\leq 2|\alpha|\tag{14} $$ (see e.g.
Lemma 4.2 in \cite{Agl-McCar})).

We shall prove the existence of a mapping $F=F_{\omega}\in\Cal
O(G_2,\Bbb D)$, where $\omega\in\partial\Bbb D$, such that $$
\lim\sb{\lambda\to
0}\frac{1}{|\lambda|}m(F\circ\phi(0),F\circ\phi(\lambda))=1,\tag{15}$$
where
$m(\lambda_1,\lambda_2):=\left|\frac{\lambda_1-\lambda_2}{1-\bar\lambda_1\lambda_2}\right|$,
$\lambda_1,\lambda_2\in\Bbb D$. The above equality implies that
$F\circ\phi$ is an automorphism of the unit disk and,
consequently, it shows that $\phi$ is a complex geodesic.

Elementary calculation show that the limit in \thetag{15} equals $$
\frac{|(1+\tau\alpha\bar\omega)^2+\tau(1-|\alpha|^2)|}{2+2\Re(\bar\omega\tau\alpha)}.
$$ The last expression is smaller than or equal to $$
\frac{|1+\tau\alpha\bar\omega|^2+1-|\alpha|^2}{2+2\Re(\bar\omega\tau\alpha)}
$$ with the equality holding iff
$\frac{(1+\tau\alpha\bar\omega)^2}{\tau}>0$. But the last expression equals
$1$. So we shall finish the proof if we show for any $\tau\in\partial\Bbb D$
satisfying \thetag{14} the existence of $\tilde\omega\in\partial\Bbb D$ such
that $\frac{(1+\tilde\omega|\alpha|)^2}{\tau}>0$. It is easy to verify that
it is sufficient to show the existence of such an $\tilde\omega$ for
$\tau\in\partial\Bbb D$ satisfying the equality $|1-\tau|=2|\alpha|$ or
$1-\cos\tau=2|\alpha|^2$ (and then without loss of generality we may assume
that $\sin\tau>0$). So let
$\tau=1-2|\alpha|^2+i\sqrt{1-|\alpha|^2}2|\alpha|$. Define
$\tilde\omega:=-|\alpha|+i\sqrt{1-|\alpha^2}$. Then we may easily verify
that $\frac{(1+\tilde\omega|\alpha|)^2}{\tau}>0$, which finishes the proof.
\qed \enddemo

\subheading{Remark} Note that a more effective description of
automorphisms $f_1$ and $f_2$ coming out in the case \thetag{ii}
is given with he help of the inequality \thetag{14}. Although the
condition \thetag{14} applies to the situation when $f_1=\id$, it
gives no restriction on generality of the characterization.

The proof of the equality $c_{G_2}=\tilde k_{G_2}$ given by
Costara in \cite{Cos~2} (unlike the one of Agler and Young in
\cite{Agl-You~2}) is purely function-theoretic. So all the results
obtained in the paper may be proven by the function theoretic
means.

Our method led us to an alternate (but certainly equivalent)
description of complex geodesics in the symmetrized bidisc to the
one of Agler, Costara, Yeh and Young. In particular, our
description suggests a possible form of $\tilde k_{G_n}$-geodesics
(or probably complex geodesics) for $G_n$, $n\geq 3$, where $G_n$
denotes the higher dimensional analogue of the symmetrized bidisc
defined as follows. $G_n:=\pi_n(\Bbb D^n)$, where $$
\pi_n(\lambda_1,\ldots,\lambda_n)=(\lambda_1+\ldots+\lambda_n,\sum_{1\leq
j<k\leq
n}\lambda_j\lambda_k,\ldots,\lambda_1\cdot\ldots\cdot\lambda_n),
$$ $\lambda_j\in\Bbb C$, $j=1,\ldots,n$. Namely, the suggested
form of the $\tilde k_{G_n}$-geodesics passing through
$(0,\ldots,0)$ (and for our convenience mapping $0$ into
$(0,\ldots,0)$) is of the following form
$$\phi(\lambda)=\pi_n\circ(B(\epsilon_0\root{n}\of{\lambda}),\ldots,B(\epsilon_{n-1}\root{n}\of{\lambda})),
\;\lambda\in\Bbb D, $$ where in the definition above $B$ is any
Blaschke product of degree less than or equal to $n$, $B(0)=0$,
$\root{n}\of{1}=\{\epsilon_0,\ldots,\epsilon_{n-1}\}$ and the root
$\root{n}\of{\lambda}$ is chosen arbitrarily. The main obstacle in
the possible proof of the fact that these functions are actually
$\tilde k_{G_n}$-geodesics (or even complex geodesics) is the lack
of good candidates of functions maximalizing the expression in the
definition of the Carath\'eodory distance, which would play the
role of the functions $F_{\omega}$ from the two dimensional case.


\Refs \widestnumber\key{Agl-Yeh-You}

\ref \key Agl-McCar \by J. Agler, J. E. McCarthy \paper Norm
preserving extensions of holomorphic functions from subvarieties
of the bidisk \jour Annals of Mathematics \vol 157 \yr 2003 \pages
289--312
\endref

\ref \key Agl-Yeh-You \by J. Agler, F. B. Yeh \paper Realization
of functions into the symmetrised bidisc \jour preprint
\endref

\ref \key Agl-You~1 \by J. Agler, N. J. Young \paper A Schwarz
Lemma for the symmetrized bidisc \jour Bull. London Math. Soc.
\vol 33 \yr 2001 \pages 195--210
\endref

\ref \key Agl-You~2 \by J. Agler, N. J. Young \paper The
hyperbolic geometry of the symmetrized bidisc \jour preprint
\endref

\ref \key Cos~1 \by C. Costara \paper The symmetrized bidisc as a
counterexample to the converse of Lempert's theorem \jour Bull.
London Math. Soc. \toappear
\endref

\ref \key Cos~2 \by C. Costara \paper Dissertation \jour
Universit\'e Laval\endref

\ref \key Edi~1 \by  A.  Edigarian  \paper  On  extremal mappings
in complex ellipsoids \jour Ann. Pol. Math. \vol 62 \yr 1995
\pages 83--96
\endref

\ref \key Edi~2 \by A. Edigarian \paper private communication (October 2003)
\endref

\ref \key Jar-Pfl \by M.Jarnicki \& P.Pflug \book Invariant
Distances and Metrics in Complex Analysis \publ Walter de Gruyter
\yr 1993 \endref

\ref \key Jar-Pfl-Zei \by M. Jarnicki, P. Pflug \& R. Zeinstra
\paper Geodesics for convex complex ellipsoids \jour Ann. d. Sc.
Norm. Sup. di Pisa \vol XX Fasc 4 \yr 1993 \pages 535--543
\endref

\ref \key Kob \by S. Kobayashi \book Hyperbolic Complex Spaces
\publ Grundlehren der mathematischen Wissenschaften, Springer, vol
318 \yr 1998
\endref

\ref \key Lem \by L.Lempert \paper La m\'etrique de Kobayashi et
la repr\'esentation des domaines sur la boule \jour Bull. Soc.
Math. France \vol 109 \yr 1981 \pages 427-479
\endref

\ref \key Pfl-Yous \by P. Pflug, E. H. Youssfi \paper Complex
geodesics of the minimal ball in $\Bbb C^n$ \jour preprint \endref

\ref \key Pfl-Zwo \by P. Pflug \& W. Zwonek \paper The Kobayashi metric for
non-convex complex ellipsoids \jour Complex Variables \vol 29 \yr 1996
\pages 59-71
\endref

\ref \key Zwo \by W. Zwonek \paper Effective  formulas for complex
geodesics in generalized pseudoellipsoids with applications \jour
Ann. Pol. Math. \vol 61 \yr 1995 \pages 261--294
\endref
\endRefs
\enddocument